# Non-conventional Estimation Theorems Concerning a Ubiquitous Bilinear Stochastic Differential System: a Control Perspective

Sandhya Rathore

*Department of Electrical Engineering, Sardar Vallabhbhai National Institute of Technology, Surat, Gujarat-395007, India. E-mail address: sandhya.rathore.svnit@gmail.com*

Shambhu Nath Sharma

*Department of Electrical Engineering, Sardar Vallabhbhai National Institute of Technology, Surat, Gujarat-395007, India. E-mail address: snsvolterra@gmail.com*

Đani Juričić

*Department of Systems and Control, Jožef Stefan Institute, Jamova cesta 39, SI-1000 Ljubljana, Slovenia E-mail address: juricic@ijs.si*

**Abstract**

Stochastic differential equations and the associated partial differential equations are the cornerstone formalism in stochastic control problems. The universality of bilinear stochastic systems can be found in autonomous systems, non-linear dynamic circuits, and mathematical finance. Consensus on the Itô versus Stratonovich suggests 'stochastic systems embedded with 'Stratonovich differential' to describe the stochastic evolution of the state vector of real physical systems'. The mathematical theory of a *scalar* time-varying bilinear Stratonovich stochastic differential equation is available in current texts. The theory of scalar Stratonovich systems was developed by deriving their closed-form solutions and then conditional moments.

Practical problems obeying *vector* time-varying bilinear Stratonovich stochastic differential equations are ubiquitous. However, their formal and systematic estimation theory is not available. In this paper, we develop the results for non-homogeneous Markov processes obeying the vector Stratonovich bilinear stochastic differential equation. Then, the theory of the paper is applied to a stochastic three-phase rectifier circuit. The stochastic evolution of the rectifier state vector is constructed by utilizing the Euler-Lagrange theory and embedding Stratonovich differential. Since the theory of the paper accounts for the multi-dimensionality as well as respects Stratonovich stochasticity of the ubiquitous bilinear system, the estimation Theorems will be remarkably useful to estimation and control of non-linear real physical systems.

Finally, this paper will be of interest to control and computational practitioners aspiring for the formal theory of bilinear stochastic systems arising from their application as well as applied mathematicians looking for applications of formal bilinear stochastic estimation theory.

## 1. Introduction

In the theory of stochastic processes, the Itô theory has found its striking applications to model the random forcing term of dynamical systems, see Kunita [1]. In the Itô theory, the term $dW_t$ is a rigorous mathematical object, where $W_t$ is the Brownian motion process. The rigorous formulation of the Brownian motion process can be traced back to pioneering works of Norbert Wiener and A N Kolmogorov, see Sharma and Gawalwad [2]. Multi-dimensional stochastic differential rules and the characteristic function have proven useful to construct the theory of stochastic differential systems from a control perspective (Karatzas and Shreve, [3], p. 384). Notably, the multi-dimensional stochastic differential rule assumes the structure of the Stochastic Differential Equation (SDE). That can be regarded as the stochastic evolution of the scalar function. In systems theory, the scalar non-linear function of the state vector has central importance, which has interpretations as the Lyapunov function, the Lagrangian, rate functions in large deviations and conditional characteristic function. The evolution of conditional characteristic function is a consequence of the action of conditional expectation operator on the multi-dimensional stochastic differential rule. In this paper, we restrict our discussions to vector stochastic differential equations. Publications on bilinear stochastic differential equations can be found in Brockett [4], Wilsky and Marcus [5] and Zhang [6]. The closed-form stationary solutions to *scalar* bilinear Stratonovich SDE as well as the Itô SDE can be found in recent papers of Terdik [7], Patil and Sharma [8]. In Patil and Sharma [8], *scalar* Stratonovich bilinear stochastic differential equations were the subject of investigations and their estimation theory was developed. Despite Stratonovich *variation* on the Itô's theme (Stroock [9]) and the consensus (Mannella and McClintock [10]), fecundity of multi-dimensional Stratonovich bilinear stochastic differential equations in appealing practical problems, publications on formal, systematic and unified estimation theory of Stratonovich stochastic differential systems are not *explicitly* available in literature yet. The practical Stratonovich bilinear systems are helicopter rotor dynamics under turbulence (Kloeden and Platen [11]), switched electrical networks under switching uncertainty. Here, we take a little pause and explain succinctly the fecundity of multi-dimensional Stratonovich stochastic differential equations by considering a PWM rectifier under the stochastic influence. The stochastic rectifier circuit is appealing and non-trivial in circuits and systems. The problem of achieving estimations of the stochastic rectifier circuit involves the Lagrangian setting (Scherpen *et al.* [12]), Stratonovich differential (Stroock [9]) and conditional characteristic function. This example suggests 'circuit theory meets systems theory'.

The intent of this paper is to develop the theory of a multi-dimensional bilinear Stratonovich stochastic differential equation and demonstrate an application of the theory to a practical problem described via the Stratonovich stochastic differential equation. The term 'Non-conventional' of the paper means 'we begin from the 'Stratonovich' setting (Stroock [9]) and then, we construct the mathematics of the scalar function $d\exp(s^T x_t)$, where the state vector $x_t$ is a non-homogeneous Markov process obeying the vector Stratonovich time-varying stochastic differential equation. Finally, we sketch two estimation Theorems for the vector Stratonovich stochastic differential equation and they are applied to a stochastic rectifier circuit. The idea of the method of the paper bears a resemblance to a non-linear filtering method of Liptser and Shirayev [13] for statistics of random processes, see FKK filtering [14] as well. In systems and control literature, control of the rectifier circuit was realized via the Lagrangian setting. That can be found in a recently published appealing paper, i.e. Scherpen *et al.* [12]. In contrast to [12], we arrive at the stochastic differential equation description for the rectifier circuit by introducing the notion of Lagrangian formulation and Stratonovich differential. Then, we sketch the estimation procedure of the rectifier stochastic differential equation. This paper can be regarded as an *extension* of a recently published paper by one of us with another co-author, i.e. Patil and Sharma [8]. In [8], the estimation theory of a *scalar* Stratonovich bilinear stochastic differential equation was developed. On the other hand, this paper develops the estimation theory of a multi-dimensional Stratonovich bilinear stochastic differential equation. Thus, our paper is different and generalization of recently published papers in systems and control literature, i.e. Patil and Sharma [8] and Scherpen *et al*. [12]. The main results of this paper are more general in contrast to [8] in the sense that the main results embed the greater class of bilinear systems.

*Notations:* Throughout the paper, the sign '∘' denotes the Stratonovich interpretation. The terms $(B_0^i(t) + \sum_\phi B_{i\phi} x_\phi) \circ dW_t$ and $(B_0^i(t) + \sum_\phi B_{i\phi} x_\phi) dW_t$ denote the Stratonovich and Itô differentials respectively. The notation $\langle \ \rangle$ denotes the conditional expectation operator. Since this paper discusses the *vector* version of the Stratonovich bilinear stochastic differential equation, 'appropriate variables' for subscripts and superscripts are chosen for the component-wise description of essential formalisms of the paper.

## 2. Bilinear Estimation Theorems

Challenges to control can be circumvented by developing greater understanding of mathematical system theory. Notably, bilinear systems found their striking application in diverse fields. David Mumford argues the dawning of the age of stochasticity. Secondly, Stratonovich stochasticity is imperative for real physical systems. That obeys consensus on the Itô versus Stratonovich differentials. First, we develop a Lemma for the vector Stratonovich bilinear stochastic differential equation. Then, the Lemma will be applied to develop two estimation Theorems for the Stratonovich bilinear system. Finally, the two Theorems will be applied to a stochastic rectifier circuit to demonstrate the usefulness of the Theorems. A Lemma, two Theorems and their application to a non-trivial stochastic rectifier circuit constitute the main results of the paper. Since this paper is intended to develop the theory of the Stratonovich stochastic differential system: vector case, we weave the advanced proofs of the Theorems of the paper in component-wise setting. Here, the 'vector' means the Stratonovich state $x_t$ (Stroock [9]) is a set of $n$ random variables, where each random variable is real-valued with its stochastic evolution. The usefulness of the main results of the paper from the control perspective is to achieve the control of the bilinear state vector by exploiting estimated bilinear states to generate control signals, where observables are not available.

**Lemma**: Consider the vector Stratonovich bilinear SDE, i.e.

$$dx_t = (A_0(t) + A_t x_t)\, dt + (B_0(t) + B_t x_t) \circ dW_t,$$

and the $i$ th component of the vector Stratonovich bilinear SDE,

$$dx_i(t) = (A_0^i(t) + \sum_\alpha A_{i\alpha}(t) x_\alpha(t))\, dt + (B_0^i(t) + \sum_\phi B_{i\phi} x_\phi) \circ dW_t$$

$$= (A_0^i(t) + \sum_\alpha A_{i\alpha}(t) x_\alpha(t) + \frac{1}{2} \sum_q (B_0^q(t) + \sum_\phi B_{q\phi}(t) x_\phi(t)) B_{iq}(t)) dt$$

$$+ (B_0^i(t) + \sum_\phi B_{i\phi} x_\phi) dW_t.$$

Suppose the scalar input Brownian motion process= $\{W_t \in R, 0 \le t < \infty\}$. Then, the $(i,j)$ component of the diffusion matrix, $(bb^T)_{ij}(x_t,t)$, of the bilinear stochastic system and its second-order partial derivative become

$$(bb^T)_{ij}(x_t,t) = B_0^i(t)B_0^j(t) + \sum_\gamma B_0^i(t)B_{j\gamma}(t)x_\gamma(t) + \sum_\phi B_0^j(t)B_{i\phi}(t)x_\phi(t)$$
$$+ \sum_{\phi,\gamma} B_{i\phi}(t)B_{j\gamma}(t)x_\phi(t)x_\gamma(t)), \tag{1a}$$

$$\frac{\partial^2}{\partial \langle x_p \rangle \partial \langle x_q \rangle}(bb^T)_{ij}(\langle x_t \rangle, t) = B_{ip}(t)B_{jq}(t) + B_{iq}(t)B_{jp}(t) \tag{1b}$$

respectively.

***Proof***: Here, we compute the diffusion coefficient matrix of a vector time-varying Stratonovich bilinear stochastic differential equation. The $i$ th component of the vector Stratonovich SDE is

$$dx_i(t) = (A_0^i(t) + \sum_\alpha A_{i\alpha}(t)x_\alpha(t)) dt + (B_0^i(t) + \sum_\phi B_{i\phi} x_\phi) \circ dW_t$$
$$= (A_0^i(t) + \sum_\alpha A_{i\alpha}(t)x_\alpha(t) + \frac{1}{2}\sum_q (B_0^q(t) + \sum_\phi B_{q\phi}(t)x_\phi(t))B_{iq}(t))dt$$
$$+ (B_0^i(t) + \sum_\phi B_{i\phi} x_\phi) dW_t. \tag{2}$$

The right-hand side of the above equation is a consequence of the Itô-Stratonovich integral. The relationship between two stochastic integrals is a consequence of the mean square convergence (Karatzas and Shreve [3], p. 156). A proof of the Stratonovich correction term $\frac{1}{2}\sum_j (B_0^j(t) + \sum_\phi B_{j\phi}(t)x_\phi(t))B_{ij}(t)$ is explained in the *appendix* of the paper. In this paper, we adopt the *notations* of Karatzas and Shreve [3] to state the stochastic differential equation. The terms $a(x_t,t)$ and $a'(x_t,t)$ denote the drift and the Stratonovich correction term of non-linear stochastic differential equations. The terms $b(t,x_t)$, $(bb^T)(t,x_t)$ have interpretations as the process noise coefficient and diffusion coefficient matrix respectively. Thus,

$$a_i(x_t,t) = A_0^i(t) + \sum_\alpha A_{i\alpha}(t)x_\alpha(t), \quad a'_i(x_t,t) = \frac{1}{2}\sum_q (B_0^q(t) + \sum_\phi B_{q\phi}(t)x_\phi(t))B_{iq}(t),$$
$$b_i(t,x_t) = (B_0^i(t) + \sum_\phi B_{i\phi} x_\phi) dW_t.$$

In the component-wise setting, the diffusion coefficient matrix component $(bb^T)_{ij}(x_t,t)$ associated with equation (2) obeys the following relationship:

$$dx_i dx_j = (bb^T)_{ij}(x_t,t)dt = ((A_0^i(t) + \sum_\alpha A_{i\alpha}(t)x_\alpha(t) + \frac{1}{2}\sum_q (B_0^q(t) + \sum_\phi B_{q\phi}(t)x_\phi(t))B_{iq}(t))dt$$
$$+ (B_0^i(t) + \sum_\phi B_{i\phi} x_\phi) dW_t) ((A_0^j(t) + \sum_\alpha A_{j\alpha}(t)x_\alpha(t)$$
$$+ \frac{1}{2}\sum_q (B_0^q(t) + \sum_\phi B_{q\phi}(t)x_\phi(t))B_{jq}(t))dt + (B_0^j(t) + \sum_\gamma B_{j\gamma} x_\gamma)dW_t)$$
$$= (B_0^i(t) + \sum_\phi B_{i\phi}(t)x_\phi(t)) (B_0^j(t) + \sum_\gamma B_{j\gamma}(t)x_\gamma(t))dt$$
$$= B_0^i(t)B_0^j(t) + \sum_\gamma B_0^i(t)B_{j\gamma}(t)x_\gamma(t) + \sum_\phi B_0^j(t)B_{i\phi}(t)x_\phi(t)$$

$$+ \sum_{\phi,\gamma} B_{i\phi}(t) B_{j\gamma}(t) x_\phi(t) x_\gamma(t)) dt.$$

After comparing the both sides of the above expression, we arrive at equation (1a). Now, we compute the second-order partial of the $(i, j)$th component of the diffusion coefficient matrix $(bb^T)(x_t, t)$ evaluated at its conditional mean $\langle x_t \rangle$, i.e.

$$\frac{\partial^2}{\partial \langle x_p \rangle \partial \langle x_q \rangle} (bb^T)_{ij}(\langle x_t \rangle, t) = \frac{\partial^2}{\partial \langle x_p \rangle \partial \langle x_q \rangle} (B_0^i(t) B_0^j(t) + \sum_\gamma B_0^i(t) B_{j\gamma}(t) \langle x_\gamma(t) \rangle$$

$$+ \sum_\phi B_0^j(t) B_{i\phi}(t) \langle x_\phi(t) \rangle + \sum_{\phi,\gamma} B_{i\phi}(t) B_{j\gamma}(t) \langle x_\phi(t) \rangle \langle x_\gamma \rangle)$$

$$= \frac{\partial^2}{\partial \langle x_p \rangle \partial \langle x_q \rangle} (\sum_{\phi,\gamma} B_{i\phi}(t) B_{j\gamma}(t) \langle x_\phi(t) \rangle \langle x_\gamma \rangle)$$

$$= \sum_{\phi,\gamma} B_{i\phi}(t) B_{j\gamma}(t) \frac{\partial^2}{\partial \langle x_p \rangle \partial \langle x_q \rangle} \langle x_\phi \rangle \langle x_\gamma \rangle$$

$$= \sum_{\phi,\gamma} B_{i\phi}(t) B_{j\gamma}(t) (\delta_{\phi p} \delta_{\gamma q} + \delta_{\gamma p} \delta_{\phi q})$$

$$= B_{ip}(t) B_{jq}(t) + B_{iq}(t) B_{jp}(t).$$

Thus, we get equation (1b). Equations (1a) and (1b) of the Lemma are useful to prove the estimation Theorems of the paper. This completes the proof.

*QED*

*Remark* 1: The proof of the Lemma is about the diffusion coefficient matrix and its double derivative coupled with the *scalar input* Brownian motion process. On the other hand, the diffusion coefficient matrix and its double derivative coupled with the vector input Brownian motion process in the component-wise description become

$$\sum_\phi B_0^{i\phi}(t) B_0^{j\phi}(t) + x_i \sum_\phi B_0^{j\phi}(t) B_\phi(t) (\delta_{jq} \delta_{ip} + \delta_{iq} \delta_{jp}) \sum_\phi B_\phi^2(t) + \sum_\phi B_0^{i\phi}(t) B_\phi(t) x_j(t) + x_i x_j \sum_\phi B_\phi^2(t),$$

$$\frac{\partial^2}{\partial \langle x_p \rangle \partial \langle x_q \rangle} (bb^T)_{ij}(\langle x_t \rangle, t) = \frac{\partial^2}{\partial \langle x_p \rangle \partial \langle x_q \rangle} (\sum_\phi B_0^{i\phi}(t) B_0^{j\phi}(t) + \langle x_i \rangle \sum_\phi B_0^{j\phi}(t) B_\phi(t)$$

$$+ \sum_\phi B_0^{i\phi}(t) B_\phi(t) \langle x_j(t) \rangle + \langle x_i \rangle \langle x_j \rangle \sum_\phi B_\phi^2(t))$$

$$= (\delta_{jq} \delta_{ip} + \delta_{iq} \delta_{jp}) \sum_\phi B_\phi^2(t)$$

respectively. The procedure is similar and we omit the intermediate steps for the vector input Brownian motion case.

*Remark* 2: A simplification of the diffusion coefficient matrix exploits the multiplication rule for the time differential and the Brownian motion differential, where the square of the time differential and the time differential-Brownian motion differential product vanish. The square of the Brownian motion differential becomes the time differential. Stroock [9] is a good source for stochastic differential multiplications, the Itô and Stratonovich processes. For the similar multiplication rule for the quantum stochastic differential equation, Parthasarathy [13] can be consulted.

**Theorem 1:** The evolution of the conditional characteristic function for the vector Stratonovich time-varying bilinear stochastic differential equation, i.e.

$$dx_i(t) = (A_0^i(t) + \sum_\alpha A_{i\alpha}(t) x_\alpha(t)) \, dt + (B_0^i(t) + \sum_\phi B_{i\phi} x_\phi) \circ dW_t$$

$$= (A_0^i(t) + \sum_\alpha A_{i\alpha}(t) x_\alpha(t) + \frac{1}{2} \sum_q (B_0^q(t) + \sum_\phi B_{q\phi}(t) \langle x_\phi(t) \rangle) B_{iq}(t)) dt$$

$$+ (B_0^i(t) + \sum_\phi B_{i\phi} x_\phi) \, dW_t,$$

becomes

$$d \langle \exp(s^T x_t) \rangle = \Big\langle \sum_p (A_0^p(t) + \sum_\alpha A_{p\alpha}(t) x_\alpha(t)$$

$$+ \frac{1}{2} \sum_q (B_0^q(t) + \sum_\phi B_{q\phi}(t) x_\phi(t)) B_{pq}(t)) s_p \exp(s^T x_t)$$

$$+ \frac{1}{2} \sum_{p,q} (B_0^p(t) B_0^q(t) + \sum_\gamma B_0^p(t) B_{q\gamma}(t) x_\gamma(t) + \sum_\phi B_0^q(t) B_{p\phi}(t) x_\phi(t)$$

$$+ \sum_{\phi,\gamma} B_{p\phi}(t) B_{q\gamma}(t) x_\phi(t) x_\gamma(t)) s_p s_q \exp(s^T x_t) \Big\rangle dt.$$

***Proof:*** The proof of the conditional characteristic function evolution hinges on the property of the expectation operator on the stochastic process $\exp(s^T x_t)$. For the notational brevity and convenience, we choose the notation $\langle \ \rangle$ for the action of the conditional expectation operator $E(.)$ on lengthier stochastic terms. The notations $\langle \exp(s^T x_t) \rangle$ and $E(\exp(s^T x_t) | x_{t_0}, t_0)$ have the same interpretation.

$$\langle d \exp(s^T x_t) \rangle = \langle d\varphi(x_t, s) \rangle = \langle \varphi(x_{t+dt}, s) - \varphi(x_t, s) \rangle = \psi(t+dt, s) - \psi(t, s).$$

In the sense of characteristic function, the variable $s$ becomes $j\omega$, the input argument of the Fourier Transform. Now,

$$\langle \exp(s^T x_t) \rangle = \psi(t, s) = \langle \varphi(s, x_t) \rangle,$$

$$d \langle \exp(s^T x_t) \rangle = \psi(t+dt, s) - \psi(t, s) = \langle d \exp(s^T x_t) \rangle. \tag{3}$$

The above relation reveals the fact that the differential and the conditional expectation operator can be interchanged. Note that for deterministic initial conditions, the operator $\langle . \rangle$ becomes the expectation operator and random initial conditions, the operator $\langle . \rangle$ becomes the conditional expectation. Furthermore, the term $d \exp(s^T x_t)$ for the vector time-varying Stratonovich bilinear stochastic differential equation becomes the following stochastic evolution:

$$d \exp(s^T x_t) = (\sum_P (A_0^p(t) + \sum_\alpha A_{p\alpha}(t) x_\alpha(t) + \frac{1}{2} \sum_q (B_0^q(t)$$

$$+ \sum_\phi B_{q\phi}(t) x_\phi(t)) B_{pq}(t)) \frac{\partial \exp(s^T x_t)}{\partial x_p}$$

$$+ \frac{1}{2} \sum_{p,q} (B_0^p(t) B_0^q(t) + \sum_\gamma B_0^p(t) B_{q\gamma}(t) x_\gamma(t) + \sum_\phi B_0^q(t) B_{p\phi}(t) x_\phi(t)$$

$$+ \sum_{\phi,\gamma} B_{p\phi}(t) B_{q\gamma}(t) x_\phi(t) x_\gamma(t)) \frac{\partial^2 \exp(s^T x_t)}{\partial x_p \partial x_q}) dt$$

$$+ \sum_p (B_0^p(t) + \sum_\phi B_{p\phi} x_\phi) \frac{\partial \exp(s^T x_t)}{\partial x_p} dW_t.$$

(4)

The above stochastic differential equation is the cornerstone formalism to weave the proof of the Theorem. Alternatively, the above can be recast as

$$\exp(s^T x_t) = \exp(s^T x_{t_0}) + \sum_P \int_{t_0}^t ((A_0^p(\tau) + \sum_\alpha A_{p\alpha}(\tau) x_\alpha(\tau) + \frac{1}{2} \sum_q (B_0^q(\tau)$$

$$+ \sum_\phi B_{q\phi}(\tau) x_\phi(\tau)) B_{pq}(\tau)) \frac{\partial \exp(s^T x_\tau)}{\partial x_p}$$

$$+ \frac{1}{2} \sum_{p,q} (B_0^p(\tau) B_0^q(\tau) + \sum_\gamma B_0^p(\tau) B_{q\gamma}(\tau) x_\gamma(\tau) + \sum_\phi B_0^q(\tau) B_{p\phi}(\tau) x_\phi(\tau)$$

$$+ \sum_{\phi,\gamma} B_{p\phi}(\tau) B_{q\gamma}(\tau) x_\phi(\tau) x_\gamma(\tau)) \frac{\partial^2 \exp(s^T x_\tau)}{\partial x_p \partial x_q}) d\tau$$

$$+ \sum_p \int_{t_0}^t (B_0^p(\tau) + \sum_\phi B_{p\phi}(\tau) x_\phi(\tau)) \frac{\partial \exp(s^T x_t)}{\partial x_p} dW_\tau.$$

Note that the above equation is a consequence of the stochastic differential rule (Karatzas and Shreve [3], p. 153) in combination with equations (1a)-(2) of the paper, where $\varphi(s, x_t) = \varphi(x_t) = \exp(s^T x_t)$ with fixed $s$. From equation (3) and equation (4), we have

$$d\langle \exp(s^T x_t) \rangle = \langle (\sum_p (A_0^p(t) + \sum_\alpha A_{p\alpha}(t) x_\alpha(t) + \frac{1}{2} \sum_q (B_0^q(t)$$

$$+ \sum_\phi B_{q\phi}(t) x_\phi(t)) B_{pq}(t)) \frac{\partial \exp(s^T x_t)}{\partial x_p}$$

$$+ \frac{1}{2} \sum_{p,q} (B_0^p(t) B_0^q(t) + \sum_\gamma B_0^p(t) B_{q\gamma}(t) x_\gamma(t) + \sum_\phi B_0^q(t) B_{p\phi}(t) x_\phi(t)$$

$$+ \sum_{\phi,\gamma} B_{p\phi}(t) B_{q\gamma}(t) x_\phi(t) x_\gamma(t)) \frac{\partial^2 \exp(s^T x_t)}{\partial x_p \partial x_q}) dt$$

$$+ \sum_p (B_0^p(t) + \sum_\phi B_{p\phi} x_\phi) \frac{\partial \exp(s^T x_t)}{\partial x_p} dW_t \rangle. \tag{5}$$

Since the conditional expectation operator is a linear operator and

$$\langle \sum_p (B_0^p(t) + \sum_\phi B_{p\phi} x_\phi) \frac{\partial \exp(s^T x_t)}{\partial x_p} dW_t \rangle = 0.$$

Equation (5) boils down to

$$d\langle \exp(s^T x_t) \rangle = \langle \sum_p (A_0^p(t) + \sum_\alpha A_{p\alpha}(t) x_\alpha(t) + \frac{1}{2} \sum_q (B_0^q(t)$$

$$+ \sum_\phi B_{q\phi}(t) x_\phi(t)) B_{pq}(t)) \frac{\partial \exp(s^T x_t)}{\partial x_p}$$

$$+ \frac{1}{2} \sum_{p,q} (B_0^p(t) B_0^q(t) + \sum_\gamma B_0^p(t) B_{q\gamma}(t) x_\gamma(t) + \sum_\phi B_0^q(t) B_{p\phi}(t) x_\phi(t)$$

$$+ \sum_{\phi,\gamma} B_{p\phi}(t) B_{q\gamma}(t) x_\phi(t) x_\gamma(t)) \frac{\partial^2 \exp(s^T x_t)}{\partial x_p \partial x_q} \rangle dt. \tag{6}$$

Alternatively,

$$d\langle \exp(s^T x_t) \rangle = \langle \sum_p (A_0^p(t) + \sum_\alpha A_{p\alpha}(t) x_\alpha(t) + \frac{1}{2} \sum_q (B_0^q(t)$$

$$+ \sum_\phi B_{q\phi}(t) x_\phi(t)) B_{pq}(t)) s_p \exp(s^T x_t)$$

$$+ \frac{1}{2} \sum_{p,q} (B_0^p(t) B_0^q(t) + \sum_\gamma B_0^p(t) B_{q\gamma}(t) x_\gamma(t) + \sum_\phi B_0^q(t) B_{p\phi}(t) x_\phi(t)$$

$$+ \sum_{\phi,\gamma} B_{p\phi}(t) B_{q\gamma}(t) x_\phi(t) x_\gamma(t)) s_p s_q \exp(s^T x_t) \rangle dt. \tag{7}$$

This completes the proof. Equation (7) is a consequence of equation (6) of the paper.

*QED*

**Theorem 2:** The conditional mean and conditional variance evolutions associated with the *vector* Stratonovich bilinear stochastic differential equation stated in equation (2) are the following:

$$d\langle x_i(t) \rangle = (A_0^i(t) + \sum_\alpha A_{i\alpha}(t) \langle x_\alpha(t) \rangle + \frac{1}{2} \sum_q (B_0^q(t) + \sum_\phi B_{q\phi}(t) \langle x_\phi(t) \rangle) B_{iq}(t)) dt, \tag{8a}$$

$$dP_{ij} = (\sum_{p} P_{ip}(A_{jp}(t) + \frac{1}{2}\sum_{\lambda} B_{\lambda p}(t)B_{j\lambda}(t)) + \sum_{p} P_{jp}(A_{ip}(t) + \frac{1}{2}\sum_{\lambda} B_{\lambda p}(t)B_{i\lambda}(t))$$
$$+ B_0^i(t)B_0^j(t) + \sum_{\gamma} B_0^i(t)B_{j\gamma}(t)\langle x_{\gamma}(t)\rangle + \sum_{\phi} B_0^j(t)B_{i\phi}(t)\langle x_{\phi}(t)\rangle$$
$$+ \sum_{\phi,\gamma} B_{i\phi}(t)B_{j\gamma}(t)\langle x_{\phi}(t)\rangle\langle x_{\gamma}(t)\rangle + \sum_{p,q} B_{ip}(t)B_{jq}(t)P_{pq})dt.$$

(8b)

*Proof:* The proof of the first part of Theorem 2 hinges on the action of conditional expectation operator on the stochastic evolution of the scalar stochastic process $\exp(s^T x_t)$, where $x(t)$ is the Stratonovich state vector. Note that equation (6) can be recast using appropriate variables, i.e.

$$d\langle \exp(s^T x_t)\rangle = \Big\langle (\sum_{p}(A_0^p(t) + \sum_{\alpha} A_{p\alpha}(t)x_{\alpha}(t) + \frac{1}{2}\sum_{q}(B_0^q(t)$$

$$+ \sum_{\phi} B_{q\phi}(t)x_{\phi}(t))B_{pq}(t))\frac{\partial \exp(s^T x_t)}{\partial x_p}$$

$$+ \frac{1}{2}\sum_{p,q}(B_0^p(t)B_0^q(t) + \sum_{\gamma} B_0^p(t)B_{q\gamma}(t)x_{\gamma}(t) + \sum_{\phi} B_0^q(t)B_{p\phi}(t)x_{\phi}(t)$$

$$+ \sum_{\phi,\gamma} B_{p\phi}(t)B_{q\gamma}(t)x_{\phi}(t)x_{\gamma}(t))\frac{\partial^2 \exp(s^T x_t)}{\partial x_p \partial x_q})dt$$

$$+ \sum_{p}(B_0^p(t) + \sum_{\phi} B_{p\phi}x_{\phi})\frac{\partial \exp(s^T x_t)}{\partial x_p}dW_t\Big\rangle.$$

$$= \Big\langle \sum_{p}(A_0^p(t) + \sum_{\alpha} A_{p\alpha}(t)x_{\alpha}(t) + \frac{1}{2}\sum_{q}(B_0^q(t)$$

$$+ \sum_{\phi} B_{q\phi}(t)x_{\phi}(t))B_{pq}(t))\frac{\partial \exp(s^T x_t)}{\partial x_p}$$

$$+ \frac{1}{2}\sum_{p,q}(B_0^p(t)B_0^q(t) + \sum_{\gamma} B_0^p(t)B_{q\gamma}(t)x_{\gamma}(t) + \sum_{\phi} B_0^q(t)B_{p\phi}(t)x_{\phi}(t)$$

$$+ \sum_{\phi,\gamma} B_{p\phi}(t)B_{q\gamma}(t)x_{\phi}(t)x_{\gamma}(t))\frac{\partial^2 \exp(s^T x_t)}{\partial x_p \partial x_q} \Big\rangle dt,$$

where $\Big\langle \sum_{p}(B_0^p(t) + \sum_{\phi} B_{p\phi}x_{\phi})\frac{\partial \exp(s^T x_t)}{\partial x_p}dW_t\Big\rangle = 0$. After combining $\exp(s^T x_t) = x_i(t)$ with the above equation, we get equation (8a), i.e.

$$d\langle x_i(t)\rangle = (A_0^i(t) + \sum_\alpha A_{i\alpha}(t)\langle x_\alpha(t)\rangle + \frac{1}{2}\sum_q (B_0^q(t) + \sum_\phi B_{q\phi}(t)\langle x_\phi(t)\rangle)B_{iq}(t))dt.$$

The above conditional mean equation is different from the conditional mean equation available in traditional estimation theory, since the above accounts for the Stratonovich stochasticity. Here, we explain $\exp(s^T x_t) = x_i(t)$. The component-wise description of the matrix-vector format is adopted to construct the theory of multi-dimensional systems. Importantly, the component-wise descriptions of the conditional mean vector and the conditional variance matrix become the scalar. The term $\exp(s^T x_t)$ is also a scalar function of the state vector $x_t$, which is a general case. As the two specific cases of the term $\exp(s^T x_t)$, we consider $\exp(s^T x_t)$ as $x_i(t)$ and $x_i(t)x_j(t)$ for the conditional mean and conditional variance evolution respectively. The method can be extended for the calculation of higher-order statistics for multi-dimensional stochastic systems as well. For the proof of the second part of Theorem 2, the conditional variance relation

$$dP_{ij} = d\langle x_i x_j\rangle - d\langle x_i\rangle_i\langle x_j\rangle = d\langle x_i x_j\rangle - \langle x_i\rangle d\langle x_j\rangle - \langle x_j\rangle d\langle x_i\rangle - d\langle x_i\rangle d\langle x_j\rangle \tag{9a}$$

holds. The term $d\langle x_i x_j\rangle$ is a consequence of the conditional characteristic function evolution of equation (6), where $\exp(s^T x_t) = x_i(t)x_j(t)$. Thus,

$$\begin{aligned}d\langle x_i x_j\rangle = (&\left\langle (A_0^i(t) + \sum_\alpha A_{i\alpha}(t)x_\alpha(t) + \frac{1}{2}\sum_q (B_0^q(t) + \sum_\phi B_{q\phi}(t)x_\phi(t))B_{iq}(t))x_j \right\rangle \\ &+ \left\langle (A_0^j(t) + \sum_\alpha A_{j\alpha}(t)x_\alpha(t) + \frac{1}{2}\sum_q (B_0^q(t) + \sum_\phi B_{q\phi}(t)x_\phi(t))B_{jq}(t))x_i \right\rangle \\ &+ \left\langle B_0^i(t)B_0^j(t) + \sum_\gamma B_0^i(t)B_{j\gamma}(t)x_\gamma(t) + \sum_\phi B_0^j(t)B_{i\phi}(t)x_\phi(t) \right\rangle \\ &+ \left\langle \sum_{\phi,\gamma} B_{i\phi}(t)B_{j\gamma}(t)x_\phi(t)x_\gamma(t) \right\rangle)dt.\end{aligned} \tag{9b}$$

Note that

$$d\langle x_i(t)\rangle = (A_0^i(t) + \sum_\alpha A_{i\alpha}(t)\langle x_\alpha(t)\rangle + \frac{1}{2}\sum_q (B_0^q(t) + \sum_\phi B_{q\phi}(t)\langle x_\phi(t)\rangle)B_{iq}(t))dt, \tag{9c}$$

$$d\langle x_j(t)\rangle = (A_0^j(t) + \sum_\alpha A_{j\alpha}(t)\langle x_\alpha(t)\rangle + \frac{1}{2}\sum_q (B_0^q(t) + \sum_\phi B_{q\phi}(t)\langle x_\phi(t)\rangle)B_{jq}(t))dt, \tag{9d}$$

$$d\langle x_i(t)\rangle d\langle x_j(t)\rangle = (A_0^i(t) + \sum_\alpha A_{i\alpha}(t)\langle x_\alpha(t)\rangle + \frac{1}{2}\sum_q (B_0^q(t) + \sum_\phi B_{q\phi}(t)\langle x_\phi(t)\rangle)B_{iq}(t))$$

$$\times (A_0^j(t) + \sum_\alpha A_{j\alpha}(t)\langle x_\alpha(t)\rangle + \frac{1}{2}\sum_q (B_0^q(t) + \sum_\phi B_{q\phi}(t)\langle x_\phi(t)\rangle)B_{jq}(t))(dt)^2 = 0. \tag{9e}$$

The term $d\langle x_i(t)\rangle d\langle x_j(t)\rangle = 0$, since the square of the time differential vanishes. After embedding the above system of equations, equations (9b)-(9e), into equation (9a), and grouping the appropriate terms, we are led to

$$dP_{ij} = (\left\langle x_i(A_0^j(t) + \sum_\alpha A_{j\alpha}(t)x_\alpha(t) + \frac{1}{2}\sum_q(B_0^q(t) + \sum_\phi B_{q\phi}(t)x_\phi(t))B_{jq}(t))\right\rangle$$

$$-\left\langle x_i\right\rangle\left\langle A_o^j(t) + \sum_\alpha A_{j\alpha}(t)x_\alpha(t) + \frac{1}{2}\sum_q(B_0^q(t) + \sum_\phi B_{q\phi}(t)x_\phi(t))B_{jq}(t)\right\rangle$$

$$+\left\langle (A_0^i(t) + \sum_\alpha A_{i\alpha}(t)x_\alpha(t) + \frac{1}{2}\sum_q(B_0^q(t) + \sum_\phi B_{q\phi}(t)x_\phi(t))B_{iq}(t))x_j\right\rangle$$

$$-\left\langle A_o^i(t) + \sum_\alpha A_{i\alpha}(t)x_\alpha(t) + \frac{1}{2}\sum_q(B_0^q(t) + \sum_\phi B_{q\phi}(t)x_\phi(t))B_{iq}(t)\right\rangle\left\langle x_j\right\rangle$$

$$+\left\langle B_0^i(t)B_0^j(t) + \sum_\gamma B_0^i(t)B_{j\gamma}(t)x_\gamma(t) + \sum_\phi B_0^j(t)B_{i\phi}(t)x_\phi(t)\right\rangle$$

$$+\left\langle \sum_{\phi,\gamma} B_{i\phi}(t)B_{j\gamma}(t)x_\phi(t)x_\gamma(t)\right\rangle)dt.$$

After combining the first two terms of the right-hand side of the above equation, we get the first term of the right-hand side of equation (8b) of Theorem 2. After combining the third and fourth terms of the above equation, we get the second term of equation (8b). Finally, a simplification of the last two terms of the above equation using the conditional expectation property, we arrive at equation (8b) of Theorem 2. The conditional mean and variance evolutions of Theorem 2 are the exact evolutions of the vector Stratonovich time-varying stochastic differential equation, since they do not account for any approximation.

*QED*

*Remark 3:* The last two terms of the right-hand side of equation (8b) of Theorem 2 are the consequence of the following relation:

$$\left\langle B_0^i B_0^j(t) + \sum_\gamma B_0^i(t)B_{j\gamma}(t)x_\gamma(t) + \sum_\phi B_0^j(t)B_{i\phi}(t)x_\phi(t)\right\rangle + \left\langle \sum_{\phi,\gamma} B_{i\phi}(t)B_{j\gamma}(t)x_\phi(t)x_\gamma(t)\right\rangle$$

$$= B_0^i B_0^j(t) + \sum_\gamma B_0^i(t)B_{j\gamma}(t)\left\langle x_\gamma(t)\right\rangle + \sum_\phi B_0^j(t)B_{i\phi}(t)\left\langle x_\phi(t)\right\rangle$$

$$+ \sum_{\phi,\gamma} B_{i\phi}(t)B_{j\gamma}(t)\left\langle x_\phi(t)\right\rangle\left\langle x_\gamma(t)\right\rangle + \sum_{p,q} B_{ip}(t)B_{jq}(t)P_{pq}.$$

Furthermore, the simplification of the last term of the right-hand side, $\left\langle \sum_{\phi,\gamma} B_{i\phi}(t)B_{j\gamma}(t)x_\phi(t)x_\gamma(t)\right\rangle$, is not straightforward, we unfold the associated intermediate steps.

$$\left\langle \sum_{\phi,\gamma} B_{i\phi}(t)B_{j\gamma}(t)x_\phi(t)x_\gamma(t)\right\rangle =$$

$$= \sum_{\phi,\gamma} B_{i\phi}(t)B_{j\gamma}(t)(\langle x_\phi(t)\rangle\langle x_\gamma(t)\rangle + \frac{1}{2}\left\langle\sum_{p,q}(x_p-\langle x_p\rangle)(x_q-\langle x_q\rangle)\frac{\partial^2\langle x_\phi\rangle\langle x_\gamma\rangle}{\partial\langle x_p\rangle\partial\langle x_q\rangle}\right\rangle)$$

$$= \sum_{\phi,\gamma} B_{i\phi}(t)B_{j\gamma}(t)(\langle x_\phi(t)\rangle\langle x_\gamma(t)\rangle + \frac{1}{2}\left\langle\sum_{p,q}P_{pq}\frac{\partial^2\langle x_\phi\rangle\langle x_\gamma\rangle}{\partial\langle x_p\rangle\partial\langle x_q\rangle}\right\rangle)$$

$$= \sum_{\phi,\gamma} B_{i\phi}(t)B_{j\gamma}(t)(\langle x_\phi(t)\rangle\langle x_\gamma(t)\rangle + \frac{1}{2}\sum_{p,q}P_{pq}(\delta_{\phi p}\delta_{\gamma q}+\delta_{\phi q}\delta_{\gamma p}))$$

$$= \sum_{\phi,\gamma} B_{i\phi}(t)B_{j\gamma}(t)(\langle x_\phi(t)\rangle\langle x_\gamma(t)\rangle + P_{\phi\gamma})$$

$$= \sum_{\phi,\gamma} B_{i\phi}(t)B_{j\gamma}(t)\langle x_\phi(t)\rangle\langle x_\gamma(t)\rangle + \sum_{\phi,\gamma} B_{i\phi}(t)B_{j\gamma}(t)P_{\phi\gamma}$$

$$= \sum_{\phi,\gamma} B_{i\phi}(t)B_{j\gamma}(t)\langle x_\phi(t)\rangle\langle x_\gamma(t)\rangle + \sum_{p,q} B_{ip}(t)B_{jq}(t)P_{pq}.$$

QED

## 3. A Stratonovich rectifier circuit: An appealing example

This section revisits very briefly the Lagrangian mechanics of a PWM rectifier by writing its Lagrangian and the Euler-Lagrange equation. The Lagrangian mechanics has found its applications to develop evolution equations of autonomous systems, fluid dynamics, astrodynamics etc. This paper demonstrates an application of the Lagrangian mechanics to the rectifier circuit. The rectifier circuits have found applications, since they offer better control over dc bus voltage, low harmonic content at the source side and high power factor, see notable works of Scherpen *et al.* [12] and references therein.

A three-phase PWM rectifier is shown in figure (1), where the rectifier is supplied with the three-phase balanced supply $e_a$, $e_b$ and $e_c$ and has the internal resistance $R_i$ and the inductance $L_i$ respectively. The three switches $S_{abc}$ assume position 1 or zero when switched using the PWM technique. The parameters $R_L$ and $C$ are the load resistance and output capacitance respectively.

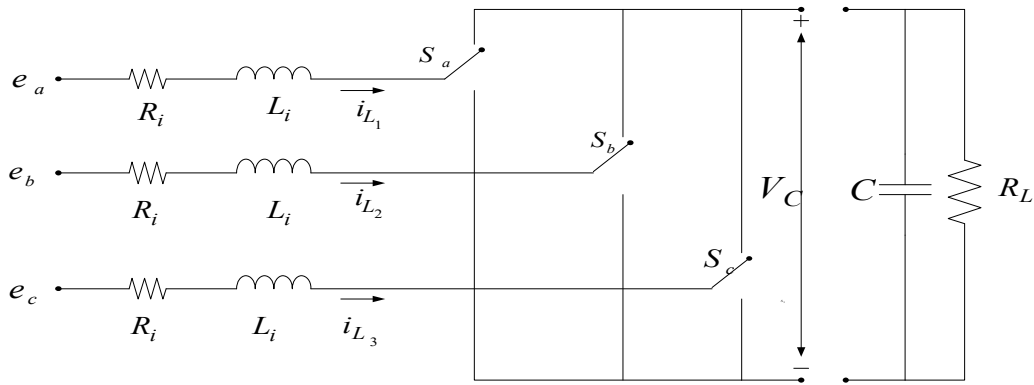

Figure 1. A schematic diagram of a three-phase rectifier

Here, we explain briefly the formulation of the Euler-Lagrange equation, the EL of the deterministic switched electrical network. Greater detail is available in a *Systems and Control* paper [12]. That can be stated as

$$\frac{d}{dt}\left(\frac{\partial \mathcal{L}}{\partial \dot{q}}(q,\dot{q})\right) - \frac{\partial \mathcal{L}}{\partial q}(q,\dot{q}) = \frac{-\partial \mathcal{D}}{\partial \dot{q}}(\dot{q}) + A(q)\lambda + \mathcal{F}_q, \quad A(q)^T \dot{q} = 0. \quad (10)$$

The Kirchhoff's current constraints are modeled by $A(q)$. Here $\dot{q}$ denotes the $n$-dimensional velocity vector and $q$ stands for the displacement, the generalized coordinates for the system. Consider $\mathcal{L}$ is the Lagrangian of the system, which represents the difference between the kinetic energy $\mathcal{T}(q,\dot{q})$ and the potential energy $\mathcal{V}(q)$ of the system. That is re-written as $\mathcal{L}(q,\dot{q}) = \mathcal{T}(q,\dot{q}) - \mathcal{V}(q)$. $\mathcal{D}(\dot{q})$ is the Rayleigh's dissipation function of the system and $\mathcal{F}_q = (\mathcal{F}_{q_1}, \mathcal{F}_{q_2}, \ldots, \mathcal{F}_{q_n})$ represents the set of generalized forcing functions associated with each generalized coordinate. Thus, space $(\mathcal{T},\mathcal{V},\mathcal{D},\mathcal{F})$ is referred to as the Lagrangian tuple of the system and modified to include the constraints introduced due to multi-loop systems, which are accounted for using the Kirchhoff's current law. Thus, the modified EL parameter space is recast as $\Sigma = (\mathcal{T},\mathcal{V},\mathcal{D},\mathcal{F},A(q))$. The EL dynamical equations are now derived for an open loop three-phase PWM rectifier under noise-free conditions. Let $\mathcal{T}$ is the switched total stored magnetic co-energy in the inductor, $\mathcal{V}$ is the switched total stored electrical energy in the capacitor. The notation $\mathcal{D}$ denotes the switched total Rayleigh dissipation function and $\mathcal{F}$ is the switched forcing function corresponding to each external source. Choosing $\dot{q}$ to denote an $n$ dimensional current vector and $q$ stands for the electrical charge. For a generalized switch position, $S_k = \{0,1\}$, we get the following for the system under considerations:

$$\mathcal{L} = \mathcal{T} - \mathcal{V} = \frac{1}{2}L_i \dot{q}_{L_a}^2 + \frac{1}{2}L_i \dot{q}_{L_b}^2 + \frac{1}{2}L_i \dot{q}_{L_c}^2 - \frac{1}{2C}q_C^2, \quad \mathcal{D} = \frac{1}{2}R_L\left(\sum_{k=a,b,c}S_k \dot{q}_{L_k} - \dot{q}_C\right)^2 + \frac{1}{2}R_i \sum_{k=a,b,c}\dot{q}_{L_k}^2,$$

$$A(q)^T \dot{q} = \dot{q}_{L_a} + \dot{q}_{L_b} + \dot{q}_{L_c} = 0, \quad e_a + e_b + e_c = 0 \quad \mathcal{F}_{q_L} = e_k, \quad \mathcal{F}_{q_C} = 0, \quad \mathcal{F}_{q_C} = 0,$$

$$k = a,b,c.$$

After combining the above system of equations with equation (10) of the paper, the rectifier current equation in $abc$ reference frame is

$$L_i \dot{i}_{L_a} = e_a - R_i i_{L_a} - (S_a - \frac{1}{3}\sum_{k=a,b,c}S_k)V_C, \quad L_i \dot{i}_{L_b} = e_b - R_i i_{L_b} - (S_b - \frac{1}{3}\sum_{k=a,b,c}S_k)V_C,$$

$$L_i \dot{i}_{L_c} = e_{3c} - R_i i_{L_c} - (S_c - \frac{1}{3}\sum_{k=a,b,c}S_k)V_C, \quad C\dot{V}_C = [S_a i_{L_a} + S_b i_{L_b} + S_c i_{L_c}] - \frac{V_C}{R_L}.$$

The above is a state space realization of the rectifier circuit. A good source of state space realizations is available in Santos *et al.* [15]. The switching signal is PWM and for sufficiently high frequency, the state space average dynamics of the PWM rectifier in $dq0$ transformation can be written as

$$L_i \dot{i}_d = e_d - i_d R_i + \omega L_i i_q - S_d V_C, \quad L_i \dot{i}_q = e_q - i_q R_i - \omega L_i i_d - S_q V_C,$$

$$C \dot{V}_{dc} = -\frac{V_C}{R_L} + (S_q i_q + S_d i_d),$$

where $(e_d \ e_q)^T = (T_{dq}^{\alpha\beta})(T_{\alpha\beta}^{abc})(e_a \ e_b \ e_c) = \left(0 \ \sqrt{\frac{2}{3}} V_m\right)^T$, The matrix $T$ refers to transformation between different frames of reference, while $V_m$ is the maximum voltage of the supply. It is assumed that the switching functions are pure sinusoids. We have

$$(S_a \ S_b \ S_c) = M\left(\sin(\omega t + \varphi_2) \ \sin(\omega t - \frac{2\pi}{3} + \varphi_2) \ \sin(\omega t + \frac{2\pi}{3} + \varphi_2)\right)^T,$$

where $(S_a \ S_b \ S_c)^T$ and $(S_d \ S_q)^T$ are switch position vectors in the $abc$ and $dq$ reference frames respectively. They are related using the following relation:

$$(S_d \ S_q)^T = (T_{dq}^{\alpha\beta})(T_{\alpha\beta}^{abc})(S_a \ S_b \ S_c) = \sqrt{\frac{2}{3}} M (\sin\varphi_2 \ \cos\varphi_2)^T.$$

For the sake of convenience, we assume $\varphi_2 = 0$ without loss of generality. The switching equivalent in the $dq0$ reference frame is explained in Han *et al.* [16], Moungkhum and Subsingha [17].

$$(S_d \ S_q)^T = \sqrt{\frac{2}{3}} M (0 \ 1)^T.$$

Note that $(i_d, i_q, V_C)^T$ is the state vector and $(e_d, e_q)^T$ represents the applied voltage vector. Thus,

$$L_i \dot{i}_d + i_d R_i - \omega L_i i_q = e_d, \quad L_i \dot{i}_q + i_q R_i + \omega L_i i_q + \sqrt{\frac{2}{3}} M V_{dc} = e_q,$$

$$C \dot{V}_C + \frac{V_C}{R_L} - \frac{3}{2}\sqrt{\frac{2}{3}} M i_q = 0. \tag{11}$$

Further, we consider the switching signal to be the gateway of stochasticity, see Sangswang and Nwankpa [18]. The noise is introduced into the system through the amplitude modulation ratio $M$. In the PWM switching, this noise effect is reflected as perturbations in the switching time. Consider $M = M(1 + \gamma \dot{W}_t)$, where M is the desired modulation index, $\gamma$ is the noise intensity and $\dot{W}_t$ is the white noise process in the context of stochastic theory, see [18]. Note that $W_t$ is the Wiener process. Other appealing system-theoretic methods, which have found applications in non-linear circuits, can be traced back to the work of Ikehara [19] using the method of Norbert Wiener. Equation (11) is further recast as

$$L_i \dot{i}_d = e_d - i_d R_i + \omega L_i i_q, \ L_i \dot{i}_q = e_q - i_q R_i - \omega L_i i_d - M\sqrt{\frac{2}{3}}(1 + \gamma \dot{W}_t) v_C,$$

$$C \dot{v}_C = -\frac{v_C}{R_L} + \sqrt{\frac{3}{2}} M(1 + \gamma \dot{W}_t) i_q.$$

We suppose a popular notation $(x_1, x_2, x_3)^T = (i_d, i_q, v_C)^T$. The above system of coupled

equations boils down to

$$\dot{x}_1 = \frac{e_d}{L_i} - \frac{R_i}{L_i}x_1 + \omega x_2, \quad \dot{x}_2 = \frac{e_q}{L_i} - \frac{R_i}{L_i}x_2 - \omega x_1 - \sqrt{\frac{2}{3}\frac{M}{L_i}}x_3 - \sqrt{\frac{2}{3}\frac{M\gamma}{L_i}}x_3\dot{W}_t,$$

$$\dot{x}_3 = -\frac{x_3}{R_L C} + \sqrt{\frac{3}{2}\frac{M}{C}}x_2 + \sqrt{\frac{3}{2}\frac{M\gamma}{C}}x_2\dot{W}_t.$$

After adopting the Stratonovich differential, we get

$$dx_t(t) = (A_0(t) + A_t x_t)dt + (B_0(t) + B_t x_t) \circ dW_t, \tag{12}$$

where

$$A_t = \begin{pmatrix} -\frac{R_i}{L_i} & \omega & 0 \\ -\omega & -\frac{R_i}{L_i} & -\sqrt{\frac{2}{3}\frac{M}{L_i}} \\ 0 & \sqrt{\frac{3}{2}\frac{M}{C}} & \frac{-1}{R_L C} \end{pmatrix}, A_0 = \begin{pmatrix} \frac{e_d}{L_i} \\ \frac{e_q}{L_i} \\ 0 \end{pmatrix}, B_0(t) = \begin{pmatrix} 0 \\ 0 \\ 0 \end{pmatrix}, B_t(t) = \begin{pmatrix} 0 & 0 & 0 \\ 0 & 0 & -\sqrt{\frac{2}{3}\frac{M}{L_i}} \\ 0 & \sqrt{\frac{3}{2}\frac{M}{C}} & 0 \end{pmatrix}.$$

For the sake of generality, we consider equation (12), a stochastic differential equation, where the parameters are time-varying. For the specific case, rectifier system parameters are not time-varying. Thus, the rectifier state vector is a Markov process, see Stroock [9]. After combining equation (6) and equation (12), we arrive at the rectifier conditional characteristic function evolution, i.e.

$$d\langle \exp(s^T x_t) \rangle = (\langle (-\frac{R_i}{L_i}x_1 + \omega x_2 + \frac{e_d}{L_i})\frac{\partial \exp(s^T x_t)}{\partial x_1} \rangle$$

$$+ \langle (-\omega x_1 - \frac{R_i}{L_i}x_2 - \sqrt{\frac{2}{3}\frac{M}{L_i}}x_3 - \frac{1}{2}\frac{M^2\gamma^2}{L_i C}x_2 + \frac{e_q}{L_i})\frac{\partial \exp(s^T x_t)}{\partial x_2} \rangle$$

$$+ \langle (-\frac{x_3}{R_L C} + \sqrt{\frac{3}{2}\frac{M}{C}}x_2 - \frac{1}{2}\frac{M^2\gamma^2}{L_i C}x_3)\frac{\partial \exp(s^T x_t)}{\partial x_3} \rangle$$

$$+ \frac{1}{2}\langle \frac{2}{3}\frac{M^2\gamma^2}{L_i^2}x_3^2\frac{\partial^2 \exp(s^T x_t)}{\partial x_2^2} + \frac{3}{2}\frac{M^2\gamma^2}{C^2}x_2^2\frac{\partial^2 \exp(s^T x_t)}{\partial x_3^2} \rangle$$

$$+ \langle -\frac{M^2\gamma^2}{L_i C}x_2 x_3 \frac{\partial^2 \exp(s^T x_t)}{\partial x_2 \partial x_3} \rangle)dt. \tag{13}$$

Suppose the probability space $(\Omega, \Im, \mu)$, where $\Omega$ is a sample space, $\Im$ is a sigma algebra and $\mu$ is a probability measure. Consider the stochastic process $x = \{x_t, \Im, 0 \le t < \infty\}$, where $\{(t, \omega)|x_t(\omega) \in B(R^d)\} \in [0, \infty) \times \Im$. Suppose the vector stochastic process $x = \{x_t, \Im, 0 \le t < \infty\}$ denotes the converter state vector. Equation (2) in combination with equation (4) describes the estimation model of the converter SDE. Then, a system of coupled condition moment equations is the following:

$$d\langle x_1\rangle = \langle -\frac{R_i}{L_i}x_1 + \omega x_2 + \frac{e_d}{L_i}\rangle dt,$$

$$d\langle x_2\rangle = \langle -\omega x_1 - (\frac{R_i}{L_i} + \frac{1}{2}\frac{M^2\gamma^2}{L_iC})x_2 - \sqrt{\frac{2}{3}}\frac{M}{L_i}x_3 + \frac{e_q}{L_i}\rangle dt,$$

$$d\langle x_3\rangle = \langle \sqrt{\frac{3}{2}}\frac{M}{C}x_2 - (\frac{1}{R_LC} + \frac{1}{2}\frac{M^2\gamma^2}{L_iC})x_3\rangle. \quad (14a)$$

$$dP_{11}(t) = 2(-\frac{R_i}{L_i}P_{11} + \omega P_{12})dt,$$

$$dP_{22}(t) = (2(-\omega P_{12} - (\frac{R_i}{L_i} + \frac{1}{2}\frac{M^2\gamma^2}{L_iC})P_{22} - \sqrt{\frac{2}{3}}\frac{M}{L_i}P_{23}) + \frac{2}{3}\frac{M^2\gamma^2}{L_i^2}\langle x_3\rangle^2 + \frac{2}{3}\frac{M^2\gamma^2}{L_i^2}P_{33})dt,$$

$$dP_{12}(t) = dP_{21} = (-\omega P_{11} - (\frac{R_i}{L_i} + \frac{1}{2}\frac{M^2\gamma^2}{L_iC})P_{12} - \sqrt{\frac{2}{3}}\frac{M}{L_i}P_{13} - \frac{R_i}{L_i}P_{12} + \omega P_{22})dt,$$

$$dP_{13}(t) = dP_{31} = (\sqrt{\frac{3}{2}}\frac{M}{C}P_{12} + (-\frac{1}{R_LC} - \frac{1}{2}\frac{M^2\gamma^2}{L_iC} - \frac{R_i}{L_i})P_{13} + \omega P_{23})dt,$$

$$dP_{23}(t) = dP_{32} = (\sqrt{\frac{3}{2}}\frac{M}{C}P_{22} - (\frac{1}{R_LC} + \frac{R_i}{L_i} + \frac{2M^2\gamma^2}{L_iC})P_{23} - \omega P_{13} - \sqrt{\frac{2}{3}}\frac{M}{L_i}P_{33} + \frac{2}{3}\frac{M^2\gamma^2}{L_i}\langle x_3\rangle^2)dt,$$

$$dP_{33}(t) = (\sqrt{\frac{3}{2}}\frac{M}{C}P_{32} - \frac{1}{R_LC}P_{33} + \frac{1}{2}\frac{M^2\gamma^2}{L_iC}P_{33} + \frac{3}{2}\frac{M^2\gamma^2}{C^2}\langle x_2\rangle^2 + \frac{3}{2}\frac{M^2\gamma^2}{C^2}P_{22})dt. \quad (14b)$$

*Remark* 4: Note that equation (14a) is a system of three coupled equations and equation (14b) is a system of six coupled equations. They constitute the rectifier coupled conditional moment equations, which are an immediate consequence of equation (12) and equations (8a)-(8b). One can arrive at the rectifier coupled conditional moments after considering equation (13) and choosing $\exp(s^T x_\tau) = x_i$ and $\exp(s^T x_\tau) = x_i x_j$ as well.

*Graphical notations*: In the simulation of the rectifier stochastic circuit, the solid line denotes the rectifier most probable trajectory, which is an immediate consequence of equation (14a). The dotted line is the actual trajectory resulting from the rectifier stochastic differential equation (12). The dash-dot is a rectifier unperturbed trajectory that obeys the rectifier ordinary differential equation (Scherpen *et al.* [12]).

## 4. Numerical Simulations

Numerical simulations are aimed to reveal the usefulness of 'bilinear estimation theory' of the paper. The bilinear Stratonovich stochastic estimation theory of the paper is formalized in two Theorems of the paper. The Theorem 2 of the paper is applied to the bilinear Stratonovich rectifier circuit. Finally, we achieve the numerical simulation of the rectifier conditional moment equations using two sets of initial conditions and system parameters. The first set of system parameters for the numerical simulation is the following:

$$f_c = 3000 \text{ Hz}, M = 0.8, \omega = 100\pi, R_i = 0.5\Omega, L_i = 1e-3 \text{ H}, V_m = 100 \text{ V},$$
$$C = 2200e\text{-}6\ \mu F, R_L = 100\ \Omega\ \gamma = 0.001.$$

The initial conditions are $i_d = 0$ Amps, $i_q = 0$ Amps $V_{dc} = 0$ Volts $P_{ii} = 0$.

The notations of the above system parameters are defined in the *preceding* section of the paper. In this paper, the bilinear Stratonovich stochastic differential equation is the subject of investigations. The bilinear Stratonovich stochastic system embeds linearity as well as bilinearity. It is important to note that Theorem (2) encompasses the qualitative characteristics of the Stratonovich bilinear stochastic system completely and the rectifier circuit state vector obeys the Stratonovich bilinear stochastic differential equation. The rectifier conditional moments are the consequence of Theorem (2) and the rectifier Stratonovich SDE. Thus, the rectifier conditional moment suggests the closeness of the estimated trajectory with the actual trajectory. Numerical simulations reveal the estimated trajectory follows the actual perturbed trajectory as well. Numerical simulations of the paper are accomplished using two different sets of data that display comparison between three trajectories.

Figures (2)-(4) suggest the rectifier most probable trajectory is bounded as well as follows the actual trajectory. The actual trajectory is attributed to the random forcing term in the Stratonovich setting. On the other hand, the conditional moment equations assume the structure of the Ordinary Differential Equation. The contribution of the random forcing term in the estimated trajectory is accounted via the Stratonovich correction term $\frac{1}{2}\sum_{q}(B_0^q(t)+\sum_{\phi}B_{q\phi}(t)x_{\phi}(t))B_{iq}(t),$ see equation (8a) of Theorem (2) of the paper. The correction term is attributed to the diffusion coefficient term of the Stratonovich bilinear stochastic differential equation. It is important to note that the conditional mean trajectory of the bilinear stochastic differential equation does not embed the random forcing term. On the other hand, the Stratonovich differential accounts for the random forcing term in the conditional mean trajectory of bilinear systems. Numerical simulations reveal 'the estimated trajectory follows the actual perturbed trajectory' as well.

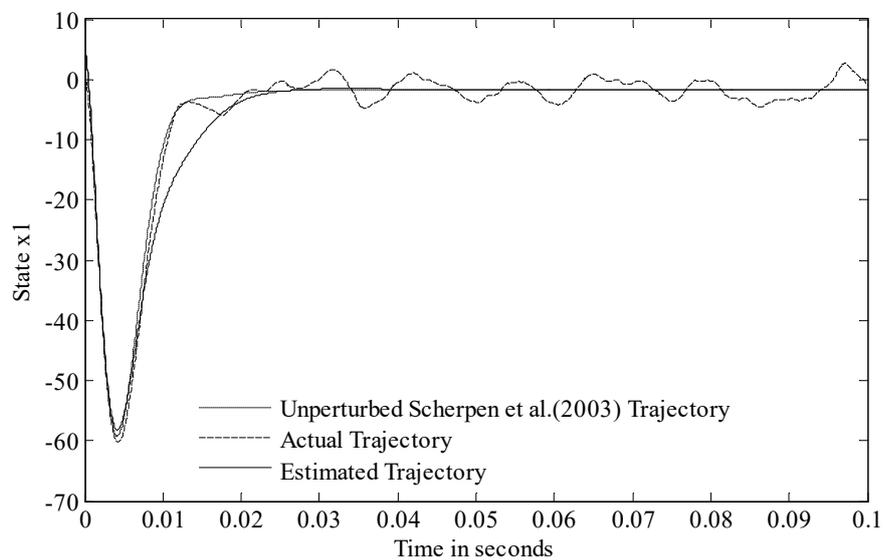

Figure 2. The state $x_1$ trajectories

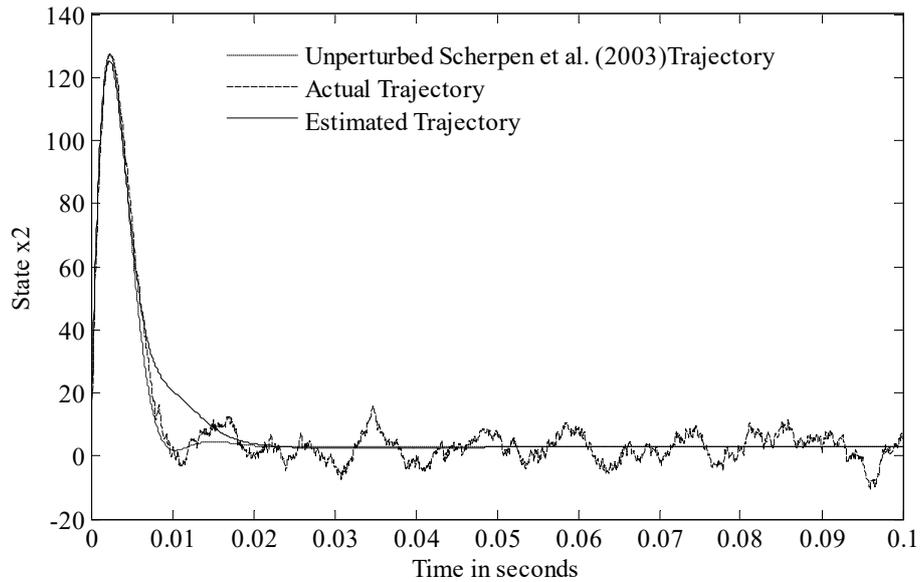

Figure 3. The state $x_2$ trajectories

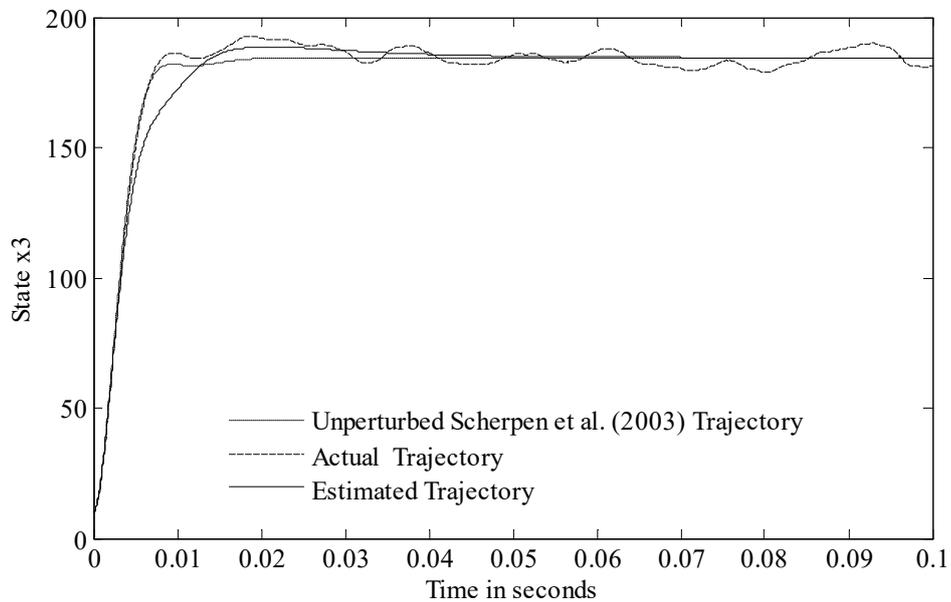

Figure 4. The state $x_3$ trajectories

Thus, the theory and numerical simulations reveal the usefulness of the main results of the paper. For the numerical simulation of the bilinear Stratonovich stochastic estimation theory of the paper, we adopt the 'estimation via numerical simulations' procedures available in Germani *et al.* [20].

In the second set of data, the efficacy of the estimator is tested by choosing a slightly higher value of $\gamma$=0.005, which goes in tandem with the system in real.

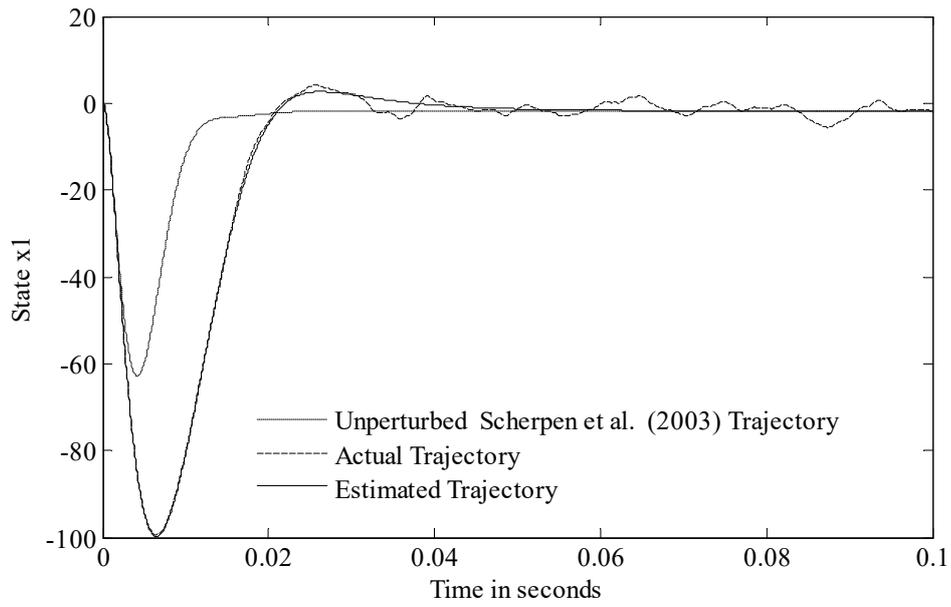

Figure 5. The state $x_1$ trajectories

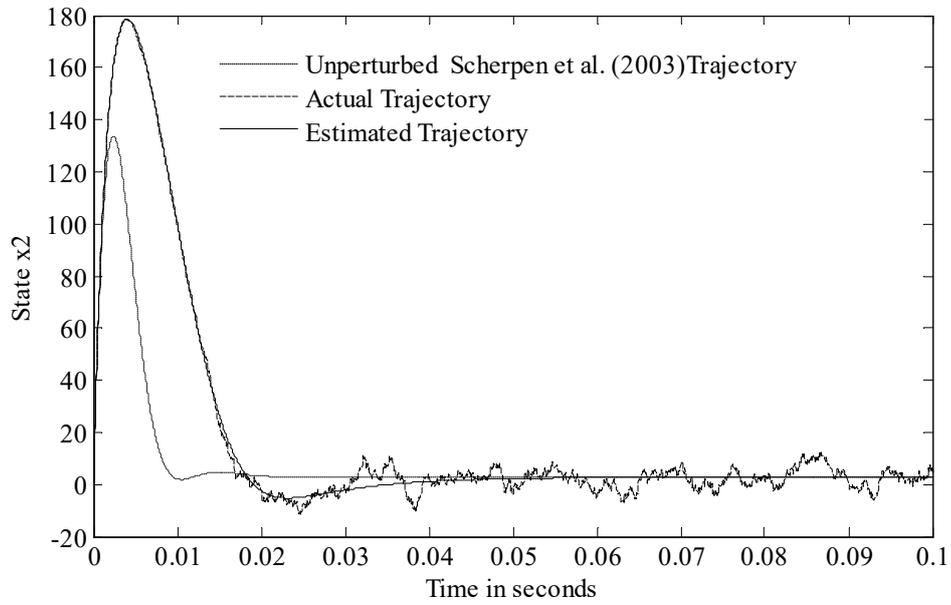

Figure 6: the state $x_2$ trajectories

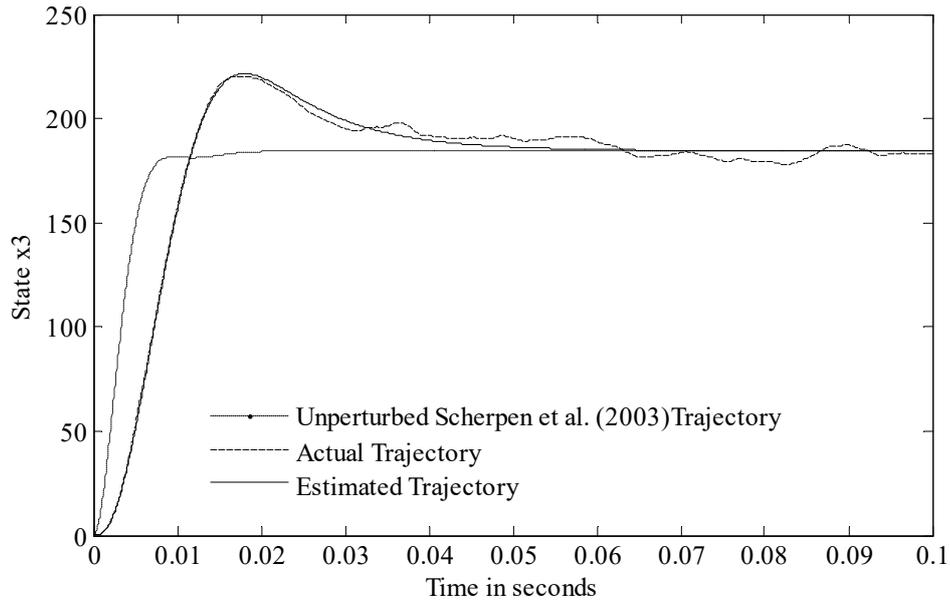

Figure 7. The state $x_3$ trajectories

Figures (5)-(7) reveal that estimated trajectories follow perturbed actual trajectories. That hold under the larger 'Stratonovich stochasticity', where unperturbed state evolutions do not agree with the most probable trajectory. Numerical simulations, displayed in figures (5)-(7), unfold the usefulness of the bilinear estimation Theorems of the paper.

**Conclusion**

The main contribution of the paper is to develop a theory of bilinear stochastic differential systems described via the *vector* time-varying Stratonovich bilinear stochastic differential equation. We have developed the theory by constructing the conditional characteristic function evolution for the Stratonovich differential. The notion of conditional characteristic function has central importance in the theory of stochastic processes. The theory of the paper is applied to a Stratonovich stochastic rectifier circuit as well. It is interesting to note that the stochastic rectifier circuit evolution *obeys* a vector Stratonovich bilinear stochastic differential equation. We have achieved the rectifier estimation procedure by developing the conditional characteristic function and conditional moment evolutions of the rectifier stochastic differential equation as well. The vector Stratonovich bilinear stochastic differential equation has an appealing structure as well as striking applications to practical problems. Another appealing application is to describe the helicopter rotor turbulence as well (Kloeden and Platen [11]).

This paper can be regarded as a *first* step towards developing a closed theory of the ubiquitous *vector* time-varying Stratonovich bilinear stochastic differential system. For the first time, the rectifier estimation, a specific case, was achieved using the method of Stratonovich differential and conditional characteristic function as well. Thus, the results of the paper are appealing. That makes a connecting thread between an unexplored connection between Stratonovich differential, bilinear systems and switched electrical networks. This paper fills a niche between applied mathematicians and control practitioners by 'sketching two formal, systematic proofs of the bilinear estimation theory and demonstrating their non-trivial application' as well. The results of *scalar* Stratonovich bilinear differential systems are available in Probability and Control Journals (Terdik [7], and Patil and Sharma, [8]). For the *vector case*, the results were not available previously despite the usefulness of the vector

Stratonovich stochastic differential equation to describe practical problems. This paper achieved that.

The control theorists and practitioners would be drawn into the fertile playground of Stratonovich stochastic systems via this paper.

**Appendix**

**Computation of Stratonovich correction terms for Bilinear Stochastic Systems**

Case I: For the sake of completeness of the paper, first we compute the Stratonovich correction term for the vector Stratonovich stochastic differential equation driven by the scalar input Brownian motion process, i.e.

$$dx_i(t) = (A_0^i(t) + \sum_\alpha A_{i\alpha}(t) x_\alpha(t)) \, dt + (B_0^i(t) + \sum_\phi B_{i\phi} x_\phi) \circ dW_t.$$

The above can be recast in a stochastic integral (Stratonovich [20], Stroock [9]), i.e.

$$x_i(t) = x_i(t_0) + \int_{t_0}^{t} (A_0^i(\tau) + \sum_\alpha A_{i\alpha}(\tau) x_\alpha(\tau)) d\tau + \int_{t_0}^{t} (B_0^i(\tau) + \sum_\phi B_{i\phi}(\tau) x_\phi(\tau)) \circ dW_\tau.$$

The first two terms of the right hand will be the same, since they have an interpretation as the Riemann integral. After applying the Stratonovich stochastic integral property to the last term of the right-hand side, we get

$$\int_{t_0}^{t} (B_0^i(\tau) + \sum_\phi B_{i\phi}(\tau) x_\phi(\tau)) \circ dW_t = LIM \sum_{k=0}^{n-1} \left( \frac{B_0^i(\tau_{K+1}) + B_0^i(\tau_k)}{2} \right.$$

$$+ \sum_\phi \frac{B_{i\phi}(\tau_{K+1}) + B_{i\phi}(\tau_k)}{2} \frac{x_0(\tau_{k+1}) + x_\phi(\tau_k)}{2} \bigg)(W_{\tau_{K+1}} - W_{\tau_k})$$

$$= LIM \sum_{k=0}^{n-1} (B_0^i(\tau_k) + \frac{B_0^i(\tau_{k+1}) - B_0^i(\tau_k)}{2}$$

$$+ \sum_\phi (B_{i\phi}(\tau_k) + \frac{B_{i\phi}(\tau_k) - B_{i\phi}(\tau_{k+1})}{2})$$

$$\times (x_\phi(\tau_k) + \frac{x_0(\tau_{k+1}) - x_\phi(\tau_k)}{2}))(W_{\tau_{K+1}} - W_{\tau_k}).$$

Note that the $LIM$ denotes the 'limit in mean', where $t_{k+1} - t_k \to 0, n \to \infty$. Thus, the above can be recast as $\int_{t_0}^{t} (B_0^i(\tau) + \sum_\phi B_{i\phi}(\tau)(x_\phi(\tau) + \frac{dx_\phi(\tau)}{2})) dW_\tau$. After taking the time derivative of $\int_{t_0}^{t} (B_0^i(\tau) + \sum_\phi B_{i\phi}(\tau)(x_\phi(\tau) + \frac{dx_\phi(\tau)}{2})) dW_\tau$, we get the stochastic term and Stratonovich correction term in the drift of the stochastic differential term. Thus, we get

$$(B_0^i(t) + \sum_\phi B_{i\phi}(t) x_\phi(t)) dW_t + \sum_\phi B_{i\phi} \frac{1}{2} dx_\phi(t) dW_t = (B_0^i(t) + \sum_\phi B_{i\phi}(t) x_\phi(t)) dW_t$$

$$+ \frac{1}{2} \sum_\phi B_{i\phi} ((A_0^\phi(t) + \sum_\alpha A_{\phi\alpha}(t) x_\alpha(t)) dt + (B_0^\phi(t) + \sum_\phi B_{\phi\gamma} x_\gamma) dW_t) dW_t$$

$$= (B_0^i(t) + \sum_\phi B_{i\phi}(t) x_\phi(t)) dW_t$$

$$+ \frac{1}{2} \sum_\phi B_{i\phi} (B_0^\phi(t) + \sum_\gamma B_{\phi\gamma} x_\gamma) dt. \qquad (A.1)$$

Note that the term $\frac{1}{2}\sum_{\phi} B_{i\phi}(B_0^{\phi}(t) + \sum_{\gamma} B_{\phi\gamma} x_{\gamma})dt$ of equation (A.1) denotes the Stratonovich correction term for the case 1. After embedding the above term into Stratonovich stochastic differential equation, we are led to the equivalent Itô differential interpretation, see equation (2) of the paper. Another appealing way of recasting the Stratonovich bilinear stochastic differential equation into the equivalent Itô differential, which is a consequence of the Stratonovich correction term in the drift term, is

$$dx_i(t) = (A_0^i(t) + \sum_{\alpha} A_{i\alpha}(t)x_{\alpha}(t) - \frac{1}{2}\sum_{q}(B_0^q(t) + \sum_{\phi} B_{q\phi}(t)x_i(t)x_i(t))B_{iq}(t))dt$$
$$+ (B_0^i(t) + \sum_{\phi} B_{i\phi} x_{\phi}) \circ dW_t$$
$$= (A_0^i(t) + \sum_{\alpha} A_{i\alpha}(t)x_{\alpha}(t))dt + (B_0^i(t) + \sum_{\phi} B_{i\phi} x_{\phi})dW_t.$$

*Case 2*: After exploiting the similar procedure for the Stratonovich bilinear stochastic differential equation for the *vector* input Brownian motion, we get the Stratonovich correction term in the drift of the Itô equivalent. For the vector input Brownian motion term, we get

$$dx_i(t) = (A_0^i(t) + \sum_{\alpha} A_{i\alpha}(t)x_{\alpha}(t))\,dt + \sum_{\phi}(B_0^{i\phi}(t) + x_i(t)B_{\phi}(t)) \circ dW_{\phi}(t)$$
$$= (A_0^i(t) + \sum_{\alpha} A_{i\alpha}(t)x_{\alpha}(t) + \frac{1}{2}\sum_{\phi}(B_0^{i\phi}(t)B_{\phi}(t) + B_{\phi}^2(t)x_i(t)))dt$$
$$+ \sum_{\phi}(B_0^{i\phi}(t) + x_i(t)B_{\phi}(t))dW_{\phi}(t). \qquad (A.2)$$

Note that the term $\frac{1}{2}\sum_{\phi}(B_0^{i\phi}(t)B_{\phi}(t) + B_{\phi}^2(t)x_i(t))$ of equation (A.2) is the Stratonovich correction term for the vector input Brownian motion case. Another appealing way of recasting the Stratonovich bilinear stochastic differential equation for the *vector* input Brownian motion into the equivalent Itô differential is

$$dx_i(t) = (A_0^i(t) + \sum_{\alpha} A_{i\alpha}(t)x_{\alpha}(t) - \frac{1}{2}\sum_{\phi}(B_0^{i\phi}(t)B_{\phi}(t) + B_{\phi}^2(t)x_i(t)))dt$$
$$+ \sum_{\phi}(B_0^{i\phi}(t) + x_i(t)B_{\phi}(t)) \circ dW_{\phi}(t)$$
$$= (A_0^i(t) + \sum_{\alpha} A_{i\alpha}(t)x_{\alpha}(t))dt + \sum_{\phi}(B_0^{i\phi}(t) + x_i(t)B_{\phi}(t))dW_{\phi}(t).$$

**Acknowledgment**: The Authors are grateful to the Researchers of the Jozef Stefan Research Institute, Ljubljana, Republic of Slovenia for providing necessary resources for completing the work embedded in this paper. The part of the work was completed during the visit of one of the Researchers to the Institute under the Joint Faculty program of the Indian National Science Academy and the Slovenian Academy of Sciences and Arts.

**References**

[1] H. Kunita, Itô's stochastic calculus: Its surprising power for applications, Stochastic Processes and their Applications. 120 (2010) 622-652.
[2] S.N. Sharma, B.G. Gawalwad, Wiener meets Kolmogorov, Norbert Wiener in the 21st Century (Thinking Machines in the Physical World), 2016 IEEE (CSS) and IEEE (SSIT) IFAC Conference, Jul 13-15 (2016) University of Melbourne, Australia.

[3] Karatzas, S.E. Shreve, Brownian Motion and Stochastic Calculus, Springer-Verlag, New York, 1991.

[4] R.W. Brockett, Modeling and Estimation with Bilinear Stochastic Systems, Defense Technical Information Center, 1979.

[5] A.S. Wilsky, S.L. Marcus, Analysis of bilinear noise models in circuits and devices, Journal of the Franklin Institute. 301 (1976) 103-122.

[6] B. Zhang, Y. Li, Model reference tracking control for discrete stochastic bilinear systems with time-varying delays, 10$^{th}$ IEEE International Conference on Control and Automation (ICCA) China, (2013) 1009-1014.

[7] G.Y. Terdik, Stationary solutions for bilinear systems with constant coefficients, in: E. Cinlar, K.L. Chug, R.K. Getoor (Eds.), Seminar on Stochastic Processes, Springer, 1989, pp. 196–206.

[8] N.S. Patil, S.N. Sharma, On the mathematical theory of a time-varying bilinear Stratonovich stochastic differential system and its application to two dynamic circuits, Transactions of the Institute of Systems, Control and Information Engineers, 27 (2014) 485-492.

[9] D. Stroock, Markov Processes from K. Itô's Perspective, Princeton University Press, 2003.

[10] R. Mannella, P.V.E. McClintock, Itô versus Stratonovich: 30 years later, Fluctuation and Noise Letters, 11 (2012) 1240010-1240019.

[11] P.E. Kloeden, E. Platen, The Numerical Solutions of Stochastic Differential Equations, Springer, New York, 1992.

[12] J.M.A Scherpen, D. Jeltsema, J.B. Klaassens, Lagrangian modeling of switched electrical networks, Systems and Control Letters. 48 (2003) 365-374.

[13] R.S. Liptser, A.N. Shiryaev, Statistics of Random Processes 1. Springer, Berlin, 1977.

[14] M. Fusisaki, G. Kallianpur, H. Kunita, Stochastic differential equations for non linear filtering problem. Osaka Journal of Mathematics, 9 (1972) 19-40.

[15] K.R. Parthasarathy, An Introduction to Quantum Stochastic Calculus, Birkhäuser, Basel, 1992.

[16] P.L.D. Santos, J.A. Ramos, T-P. Perdico´ulis, J.L. Martins de Carvalho, A note on the state-space realizations equivalence. (2011) Retrieved from : arXiv:1111.1927v1

[17] N. Moungkhum, W. Subsingha, Voltage control by DQ frame technique of SVPWM AC-DC Converter, Energy Procedia. 34 (2013) 341-350.

[18] Sangswang, C.O. Nwankpa, Effects of switching-time uncertainties on pulse width-modulated power converters: modeling and analysis, IEEE Transactions on Circuits and Systems.-I: Fundamental Theory and Applications. 50 (2003) 1006-1012.

[19] S. Ikehara, A method of Wiener in nonlinear circuits. Technical Report 217, Research Lab of Massachusetts, Cambridge, 1951.

[20] Germani, C. Manes, P. Palumbo, Linear filtering for bilinear stochastic differential systems with unknown inputs, IEEE Transactions on Automatic Control. 17 (2002) 1726-1730.

[21] R. L. Stratonovich, A new representation for stochastic integrals and equations, SIAM Journal of Control. 4 (1966) 362–371.